\documentclass{article}
\usepackage{graphicx}
\title{
Application of the Wavelet Transform with a Piecewise Linear Basis to the
Evaluation of the Hankel Transform}
\author{P. S. Zykov and E. B. Postnikov}

\date{Kursk State Pedagogical University,\\ ul. Radishcheva 33,
Kursk, 305000 Russia \\
e-mail: postnicov@mail.ru}
\begin{document}
\maketitle
\begin{abstract}
A method for computing the Hankel transform is proposed whereby the
letter is reduced to a sum by representing the integrand as a smooth function
times a Bessel function. The smooth function is replaced by its wavelet
decomposition with a basis such its scalar product with the Bessel function is
calculated analytically. The result is represented as a series, with the
coefficients strongly depending on the local behavior of the function being
transformed. The application of the method is demonstrated by an example
illustrated with plots.
\end{abstract}

The Hankel transform has found a wide variety of applications related to
problems in mathematical physics that possess an axial symmetry (see [1]). The
transforms of the zeroth ($n = 0$) and first ($n = 1$) kind,
\begin{equation}
\label{one}
\begin{array}{c}
F_n(p)=\int\limits^{\infty}_{0}f(r)J_n(pr)rdr,\\
f(r)=\int\limits^{\infty}_{0}F_n(r)J_n(pr)pdp.
\end{array}
\end{equation}
are used most frequently. This is explained by the fact that, in Heaviside's
calculus, the Laplace operator as applied to a function is equivalent to
multiplication in the Laplace transform space:
$$
\left(\frac{d^2}{dr^2}+\frac{1}{r}\frac{d}{dr}\right)f(r)
\longleftrightarrow -p^2F_0(p), \qquad
\left(\frac{d^2}{dr^2}+\frac{1}{r}\frac{d}{dr}-\frac{1}{r^2}\right)f(r)
\longleftrightarrow -p^2F_1(p),
$$

In practice, however, this method encounters a number of difficulties. One of
them is that integrals of type (1) cannot frequently be evaluated analytically
in finite form. Moreover, on computer, the function being transformed is
specified as a set of numbers. Thus, an important task is to develop numerical
methods for computing the Hankel transform with allowance made for the fact that
$J_n(pr)$ is a rapidly oscillating function for large values of $p$ (or $r$).
Two approaches to the solution of this problem are available in the literature.
The first is the fast Hankel transform, which was originally proposed in [2] and
has been improved up to the present (e.g., see Appendix B in [3] and the
references therein). This approach involves making a substitution and scaling,
which reduce the problem to logarithmic coordinates and the fast Fourier
transform. This method applies to tabulated functions, but it involves
conventional errors arising when a nonperiodic function is replaced by its
periodic extension. Additionally, the method is sensitive to the smoothness of
functions in the space of logarithmic coordinates. An efficient algorithm (also
based on the fast Fourier transform) was proposed for commuting the discrete
Hankel (Fourier--Bessel) transform in [4]. It was noted in [4] that the
algorithm is applicable in the case of an integral transform approximated by a
sum. Another method (see [5]) involves the representation of the integrand as
the product of a smooth function approximated by a set of polynomials and an
oscillating function such that the integral of its product with the polynomial
can be computed analytically. The accuracy of this approach is limited primarily
by the accuracy of the polynomial approximation of the function to be
transformed.

The basic idea behind our method is similar to that in [5]. However, instead of
a single approximation of the smooth factor, we apply its multiscale
decomposition with respect to a function basis that takes into account the local
properties of the function being transformed. Moreover, the scalar product of
basis functions must be calculated analytically.

These conditions are satisfied by some spline wavelet decompositions, which have
been intensively studied in recent years (e.g., see [6--8]). Let $g(r)$ be a
function in $L^2(R)$. Then it can be represented as a series
\begin{equation}
g(r)=\sum_{j\in Z}\sum_{k \in Z}d_{jk}\psi_{jk}(r),
\end{equation}
where $\psi_{jk}(r)=2^{j/2}\psi{2^jr-k}$ are basis functions with step $k$ and
 scale $j$ that satisfy the
conditions of spatial localization and self-similarity under scale changes. The
latter means that all $\psi_{jk}$ are formed by contracting and translating a
single function $\psi$.

Since the direct and inverse Hankel transforms (1) are symmetric, we restrict
our consideration to the direct transform and set $g(r) = f(r)r$ for simplicity.
Substituting (2) into the integral gives
\begin{equation}
F_n(p)=\sum_{j\in Z}\sum_{k \in Z}d_{jk}\int\limits^{\infty}_0
\psi_{jk}(r)J_n{pr}dr.
\end{equation}
Consider piecewise linear wavelets. Suppose that the basis wavelet is the
piecewise linear function
$$
\psi(r)=\sum_m\frac{r-mh}{h}b_{m+1}-\frac{r-(m+1)h}{-h}b_{m},
$$
where $b_m$ are the values at the junction points and $h$ is the minimum step.

The limits of summation are determined by the type of the wavelet chosen. In
particular, as a basis, one can use the Battle--Lemarie $\psi^{B,2}$
 and Str\"omberg
$\psi^{St,2}$ orthogonal wavelets [6], which are an infinite set of linear segments, or
one can apply semiorthogonal or biorthogonal bases with a compact support
consisting of $M$ segments. Examples of the latter are B-spline semiorthogonal
wavelets [8], biorthogonal wavelets obtained by applying the lifting scheme [9],
and BlaC-wavelets (Blending of Linear and Constant) [10], which combine a
piecewise linear and a piecewise constant function. Note that the Battle--
Lemarie and Str\"omberg bases rapidly decay at infinity, even though they have
no compact support. Thus, they can be truncated to $M$ segments in practical
computations.
Applying the integration formulas of [11], we find that the Hankel transform is
equivalent to a sum:
\begin{equation}
\begin{array}{l}
 F_0 (p) = \frac{1}{p}\sum\limits_{j \in Z} {\sum\limits_{k \in Z} {2^{j/2} } } d_{jk} \sum\limits_m {\left( {(b_{k + 1}  - b_k )\left[ {(k + m + 1)J_1 (2^{ - j} (k + m + 1)hp) - } \right.} \right.}  \\
 \left. { - (k + m)J_1 (2^{ - j} (k + m)ph)} \right] + \left[ {(k + m + 1)b_k  - (k + m)b_{k + 1} } \right] \times  \\
 \left\{ {(k + m + 1)J_0 (2^{ - j} (k + m + 1)hp)} \right. - (k + m)J_0 (2^{ - j} (k + m)ph) +  \\
 \left. {\left. {\frac{\pi }{2}\left[ {(k + m + 1)D(2^{ - j} (k + m + 1)hp) - (k + m)D(2^{ - j} (k + m)ph)} \right]} \right\}} \right) \\
 \end{array}
\end{equation}
\begin{equation}
\begin{array}{l}
 F_1 (p) = \frac{1}{p}\sum\limits_{j \in Z} {\sum\limits_{k \in Z} {2^{j/2} } } d_{jk} \sum\limits_m {\left( {\frac{{\pi (b_{k + 1}  - b_k )}}{2}\left[ {(k + m + 1)D(2^{ - j} (k + m + 1)hp) - } \right.} \right.}  \\
 \left. { - (k + m)D(2^{ - j} (k + m)ph)} \right] + \left[ {(k + m + 1)b_k  - (k + m)b_{k + 1} } \right] \times  \\
 \left. {\left\{ {(k + m + 1)J_0 (2^{ - j} (k + m)hp)} \right.\left. { - (k + m)J_0 (2^{ - j} (k + m + 1)ph)} \right\}} \right) \\
 \end{array}
\end{equation}
where $D(x) = H_0(x)J_1(x)-H_1(x)J_0(x)$ and $H_{0,1}$ are Struve functions [11].

Thus, if the function being transformed (which is completely characterized by
its coefficients [6--8]) has a finite-form definite integral of its product with
a linear function, then series (4) and (5) are an exact representation of the
Hankel transform as a sum. Note that for sufficiently smooth functions strongly
varying only in small ranges of their arguments, most of the $d$-coefficients in
(4) and (5) are small for high levels of resolution $j$. This makes it possible
to derive an approximate numerical Hankel transform from multiscale analysis, in
which case compactly supported bases have a significant advantage, because they
avoid additional errors associated with the finite length of the function ranges
used in numerical computations.
Suppose that the function to be transformed is specified as a set of $N$ points,
which are interpreted as its projection on an interval onto the space of scaling
functions of resolution $J\le\log_2N$.
Then (2) is replaced by the series
$$
g(r)=\sum_{j=0}^{J}\sum{k}d_{jk}\psi_{jk}(r)+
\sum_k s_{Jk}\varphi{Jk}(r),
$$
where $\varphi_{0k}$ is a scaling function of a coarser resolution.
In the computations, one can use only those of the coefficients $d_{jk}$
 that are
greater than a prescribed small number $\varepsilon$. The accuracy of the
approximation of $g$ by a truncated function $g_{\varepsilon}$
 [7] is defined as
  $\Delta=\left(\|g-g_{\varepsilon}\|\right)^{1/2}=
  \left(\left\|\sum d_{jk}^{*}\psi_{jk}\right\|^2\right)^2$ .
Taking into account the orthogonality of the basis wavelets, we obtain
$\Delta \le\varepsilon n_0^{1/2}$, where $n_0$ is the number of
discarded coefficients $d_{jk}^{*}\le \varepsilon$. The resulting
set contains wide bands of zeros in the domain where the function is smooth and
an increased number of coefficients in the domain where it undergoes abrupt
changes. Since the basis wavelet functions are narrowly localized, this implies
a high degree of adaptivity in choosing an interpolation grid. Series (4) and
(5) preserve their form when the limits in $j$ and $k$ are changed or terms are
added that correspond to scalar products with scaling functions representable in
a similar form.

Thus, the integral transforms have been reduced to sums. For small arguments of
the Bessel function, they can be computed by conventional codes and, for large
arguments, by applying algorithms developed, for example, in [4].

Numerical computations in the case where $g(r)$ is defined analytically are
interesting on their own. For this purpose, we can use an inverse version of the
lifting scheme [9] for piecewise linear basis functions.
Let a function be defined on the interval $[0,\quad 1]$ (this can always be achieved by
suitable scaling). This interval is bisected (the crudest approximation) and the
lifting scheme is applied to a set of three function values (at the endpoints
and the middle of the interval) to obtain the first iterates
$$
\gamma_{00}=g\left(\frac{1}{2}\right)-
\frac{1}{2}\sum_{k=0}^{1}g(k), \qquad \lambda_k=g(k)+\frac{1}{4}\gamma_{00},
$$
$$
V_0(x)=\sum_{k=0}^{1}\lambda_k\varphi_{0k}(x),
\qquad
W_0(x)=\gamma_{00}\psi_{00}(x),
\qquad
g(x)\approx V_0(x)+W_0(x).
$$

If the detailed coefficient $\gamma_{00}$ is not small, then the
set is supplemented with the values of $g(r)$ calculated at the
middle points of each half-interval, and the next iterates are
computed as follows. Calculate the detailed coefficient of coarser
level $J$: $$ \gamma_{jk}=g\left[2^{-(J+k)}(2k+1)\right]-
V_{0}\left[2^{-(J+k)}(2k+1)\right]-
W_{J-1}\left[2^{-(J+k)}(2k+1)\right]. $$

By using this relation, the detailed coefficients of the preceding
levels and the approximation coefficients are updated to give $$
\gamma_{jk}= \gamma_{jk}-W_J\left(2^{-(J+1)}(2k+1)\right), \quad,
 j = \overline {0,J - 1},
$$ $$ \gamma_{jk}= \gamma_{jk}-W_J\left(2^{-(J+1)}(2k+1)\right).$$
Here $W_J$ means all corresponding levels of the wavelet transform
of difference between $g(r)$ and the transform of level $J-1$.

 The current representation of the function being decomposed is
formed by
\[
\sum\limits_{k = 0}^1 {\lambda _k \varphi _{0k} (x),\quad W_j (x)
= } \sum\limits_{j = 0}^J {\sum\limits_{k = 0}^{2^j  - 1} {\gamma
_{jk} \psi _{jk} (x),\quad g(x) \approx V_0 (x) + \sum\limits_{j =
0}^J {W_j (x);} } }
\]

 Here, the scaling function and the basis wavelet are given by $$
\varphi_{0k}(x)=\Lambda(x-k), \quad
\phi(x)=\Lambda(2x-1)-\frac{1}{4}\Lambda(x)-\frac{1}{4}\Lambda(x-1),
$$ where $\Lambda(x) = \max \{0, 1-|x|\}$.

As an example, we consider a function that has an exact first-kind Hankel
transform:
$$
f(r)=re^{-ar^2}, \quad F_1(p)=\frac{p}{(2a^2)^2}\exp\left(-\frac{p^2}{4a^2}\right).
$$

Let $a = 0.4$. Then the basic nonzero part of the function is concentrated on the
interval $[0,\quad 8]$. By making the substitution $r^* =r/8$, this interval is
transformed into a unit interval, and the above procedures for wavelet
decomposition with updating and the Hankel transform (5) are applied to
 $g(r^*)=f(r^*)r^*$.
 \begin{figure}
\label{fig1}
\includegraphics[width=\textwidth]{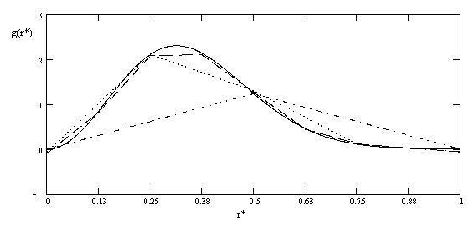}
\caption{}
\includegraphics[width=\textwidth]{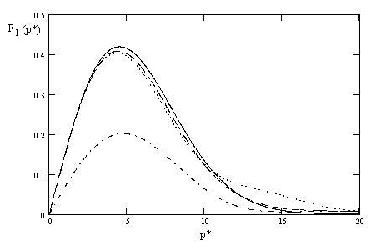}
\caption{}
 \end{figure}

  \begin{figure}
\label{fig1}
\includegraphics[width=\textwidth]{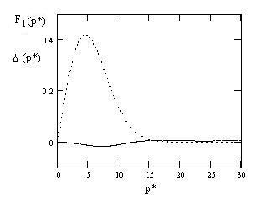}
\caption{}
 \end{figure}
Figure 1 shows the plot of $g(r^*)$ (solid curve) and its wavelet representations at
the zeroth (dot-and-dash), first (dotted), and second (dashed) levels. The
corresponding curves in Fig. 2 plot the exact and approximate Hankel transforms
in these cases. It can be seen that even this small number of levels provides an
adequate pattern (except for the zeroth approximation, which gives a rather
crude representation of the original function because of using the function
value at the middle point of the interval). The next level provides a good
approximation for most of the curve, except for large values of $p^*$. This domain
is corrected by the second approximation. Figure 3 shows the absolute error for
this level, $\Delta(p^*)$, which is equal to the difference between the
approximate (dotted) and exact (solid) solutions against the background of the
exact transform. It can be seen that the absolute error is fairly small in the
entire domain shown in the figure and does not increase sharply with $p^*$.
The computed solution is as accurate as that obtained by linear approximation at
the same nodes. However, our method provides rough estimation over the entire
domain of the function, with subsequent refinement over intervals of interest.
An improved accuracy in this approach can be achieved not only by increasing the
number of decomposition levels taken into account but also by using a basis with
unequally spaced nodes adapted to the behavior of the function to be
transformed.


\begin{thebibliography}{99}
\bibitem{1} C. J. Tranter, Integral Transforms in Mathematical Physics (Methuen,
London, 1951).
\bibitem{2} A. E. Siegman, Optics. Lett., No. 1, 13 (1977).
\bibitem{3} A. J. S. Hamilton, Mon. Not. R. Astron. Soc. 312, 257 (2000).
\bibitem{4} Ya. M. Zhileikin and A. B. Kukarkin, Zh. Vychisl. Mat. Mat. Fiz. 35, 1128
(1995).
\bibitem{5} R. Barakat and E. Parshall, Appl. Math. Lett., No. 5, 21 (1996).
\bibitem{6} S. B. Stechkin and I. Ya. Novikov, Usp Mat. Nauk 53, 53 (1998).
\bibitem{7} I. N. Dremin, O. V. Ivanov, and V. A. Nechitailo, Usp. Fiz. Nauk 171, 465
(2001).
\bibitem{8} C. K. Chui, An Introduction to Wavelets (Academic, Boston 1992; Mir,
Moscow, 2001).
\bibitem{9} W. Sweldens, in "Proceedings of Wavelet Applications in Signal and Image
Processing III" SPIE Proceedings Series 2569 (SPIE, Bellingham, 1995), p. 68.
\bibitem{10} G.-P. Bonneau, S. Hahmann, and G. M. Nielson, in "Proceedings of VIS'96,"
(ACM, New York, 1996), p. 43.
\bibitem{11} Handbook of Mathematical Functions, with Formulas, Graphs, and Mathematical
Tables, Ed. by M. Abramowitz and I. A. Stegun (Dover, New York, 1972).
\end{thebibliography}
\end{document}